\documentclass{bmcart}

\usepackage[utf8]{inputenc} 

\usepackage{amsmath,amssymb}
\makeatletter
\def\munderbar#1{\underline{\sbox\tw@{$#1$}\dp\tw@\z@\box\tw@}}
\makeatother

\newtheorem{theorem}{Theorem}

\usepackage{graphicx}  

\startlocaldefs
\endlocaldefs

\begin{document}

\begin{frontmatter}

\begin{fmbox}
\dochead{Research}


\title{Convergence of an Explicit Iterative Leap-frog Discontinuous Galerkin Method for Time-domain Maxwell's Equations in Anisotropic Materials}


\author[
   addressref={aff1},                   
   corref={aff1},                       
   email={alma@mat.uc.pt}   
]{\inits{AA}\fnm{Ad\'erito} \snm{Ara\'ujo}}
\author[
   addressref={aff1},
   email={silvia@mat.uc.pt}
]{\inits{SB}\fnm{S\'\i lvia} \snm{Barbeiro}}
\author[
   addressref={aff2},
   email={khaksar.maryam@gmail.com}
]{\inits{MKG}\fnm{Maryam} \snm{Khaksar Ghalati}}

\address[id=aff1]{
  \orgname{CMUC, Department of Mathematics, University of Coimbra}, 
  \street{Apartado 3008, EC Santa Cruz},                     %
  \city{3001-501 Coimbra},                              
  \cny{Portugal}                                    
}
\address[id=aff2]{%
  \orgname{Institute for Systems and Robotics  / LARSyS,  Instituto Superior T\'ecnico,
Universidade de Lisboa},
  \city{Lisbon},
  \cny{Portugal}
}


\begin{artnotes}
\note[id=n1]{Equal contributor} 
\end{artnotes}

\end{fmbox}


\begin{abstractbox}

\begin{abstract} 

We propose an explicit iterative leap-frog discontinuous Galerkin method for time-domain Maxwell's equations in anisotropic materials and derive its convergence properties. The {\it a priori} error estimates are illustrated  by numerical means in some experiments. 
 Motivated by a real application which encompasses modeling electromagnetic wave's propagation through the eye's structures, we simulate our model in a 2D domain aiming to represent a simple example of light scattering in the outer nuclear layer of the retina.
\end{abstract}


\begin{keyword}
\kwd{Maxwell's equations}
\kwd{explicit iterative leap-frog discontinuous Galerkin method}
\kwd{convergence}
\kwd{light scattering}
\end{keyword}


\end{abstractbox}
%

\end{frontmatter}



\section{Introduction}
\label{sec:1}

\begin{figure}[h!]
\includegraphics[scale=.5]{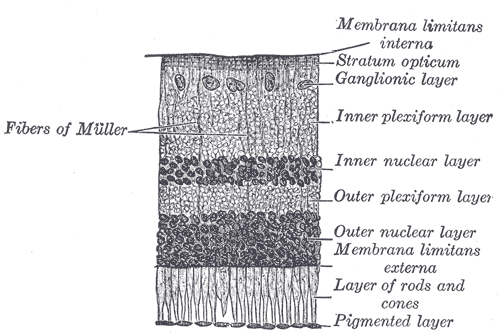}
  \caption{\csentence{Section of the retina.} Henry Gray. Anatomy of the Human Body. Philadelphia: Lea \& Febiger, 1918. (in public domain at Bartleby.com)}
\label{fig:1}      
\end{figure}

The human retina is a complex structure in the eye that is responsible for the sense of vision. It is part of the central nervous system and it is composed by several layers, namely  the outer nuclear layer that comprises the cells bodies of light sensitive photoreceptors cells, rods and cones (see Fig. \ref{fig:1}) \cite{junqueira2013}.  For many diseases that affect the eye, the diagnosis is not straightforward. The sensitivity of this structure makes medical analysis particularly  complicated. Most of the diagnoses are made either by direct observation, with the possible injection of dyes, to enhance certain parts of the organ, or by numbing the eye and directly measuring its inner pressure or thickness. There are a number of eye-related pathologies that can be identified by the detailed analysis of the retinal layers \cite{serranho2012}.  

Optical Coherence Tomography (OCT) is an increasingly popular noninvasive  technique that has been successfully used as a diagnostic tool in ophthalmology in the past decades. This method allows the assessment of the human retina {\it in vivo} and has been shown to provide functional information. By analysing data acquired through OCT, several retinal pathologies, such as diabetic retinopathy, or macular edema, can be detected in their early stages, before noticeable morphologic alterations on the retina  \cite{serranho2012}.  As OCT standard techniques only provide structural information \cite{schliesser2014}, it is necessary to expand OCT data analysis to account for both structural and functional information. OCT provides  also the possibility of evaluating different elements in measuring the retinal nerve fiber layer (RNFL), namely the tendency of RNFL thinning in glaucoma and other diseases that involve optic nerve atrophy.  Waveguides with induced anisotropy may worth to be considered for modeling  biological waveguides \cite{limeres2003}.  

Maxwell's equations are  a fundamental set of partial differential equations which describe electromagnetic wave interactions with materials.  The electromagnetic fields in space are classically described by two
field vectors, $E$ and $H$, respectively electric field and magnetic field. Here we shall consider the time domain Maxwell's equations in the transverse electric (TE) mode, as in \cite{konig2010}, where  the only non-vanishing components of the electromagnetic fields are $E_x$, $E_y$ and $H_z$. Using the following notation for the   vector and scalar curl operators
$$ \nabla \times H = \left(\frac{\partial H_z}{\partial y},-\frac{\partial H_z}{\partial x}\right)^T, \quad \mbox{curl } E=\frac{\partial E_y}{\partial x}-\frac{\partial E_x}{\partial y},$$
and assuming no conductivity effects,  the equations in the non-dimensional form are
\begin{equation}
  \epsilon\frac{\partial E}{\partial t}  = \nabla \times H, \qquad 
  \mu\frac{\partial H}{\partial t} = -\mbox{curl } E, \qquad \textrm{in } \Omega\times(0,T_f],\label{MaxwellTE}
  \end{equation}
where $E=(E_x,E_y)$ represents the electric field components and $H=(H_z)$ represents the magnetic field component.  These equations are set and solved on a bounded polygonal domain $\Omega \subset \mathbb{R}^2.$  
The electric permittivity of the medium,   $\epsilon$, and  the magnetic permeability of the medium, $\mu$, are varying  in space, being $\mu$ a scalar function and
$\epsilon$ an anisotropic tensor
\begin{equation}\epsilon=\begin{pmatrix} \epsilon_{xx} & \epsilon_{xy} \\ \epsilon_{yx} & \epsilon_{yy} \end{pmatrix}.\label{tensor}\end{equation}
We assume that electric permittivity tensor $\epsilon$ is symmetric and uniformly positive definite for almost every $(x,y) \in \Omega$, and  that it is uniformly bounded with a strictly positive lower bound, {\it i.e.}, there are constants  $\underline\epsilon>0$ and $\overline \epsilon> 0$ such that, for almost every $(x,y)\in \Omega$,
$\munderbar\epsilon |\xi|^2\le\xi^T\epsilon(x,y)\xi\le \overline\epsilon |\xi|^2,  \forall \xi\in\mathbb{R}^2.$ 
We also assume that there are constants $\underline \mu>0$ and $\overline \mu> 0$  such that, for almost every $(x,y) \in \Omega$, 
$\underline\mu \le\mu(x,y) \le \overline\mu.$

Equations (\ref{MaxwellTE}) must be complemented by initial conditions
$E(x,y,0) = E_0(x,y)$ and $H(x,y,0)= H_0(x,y)$, $(x,y)\in \Omega,$
and by proper boundary conditions. Motivated by our application of interest, here we consider  
 absorbing boundary conditions which mimic an open space by absorbing the incident radiation in the truncated computational domain. The first order Silver-M\"{u}ller absorbing boundary conditions (SM-ABC) are defined by
\begin{equation}\label{SM}
  n\times E = c\mu n \times (H\times n)\qquad \textrm {on } \partial\Omega,
\end{equation}
where $n=(n_x,n_y)^T$ is the unit outward normal vector to the boundary and $c$ is the speed with which a wave travels along the direction of the unit normal, defined, using the  effective permittivity $\epsilon_{eff}=\det(\epsilon)/(n^T \epsilon n)$ (see \cite{konig2010}), by  $c =1/\sqrt{\mu \epsilon_{eff}}$. 

The attention to the development of high-order accurate methods for solving time-domain
Maxwell's equations in complex geometries brings to the use of  discontinuous Galerkin (DG) methods \cite{hesthaven2008}. The one-step explicit time integration methods, like  leap-frog schemes, are
computationally efficient per update cycle and easy to implement. The leap-frog DG method in anisotropic materials which is discussed in \cite{fezoui2005} leads to a locally implicit method for the case of SM-ABC. 
In \cite{ABG2017} a fully explicit in time  leap-frog DG  method is investigated in the same framework. The error estimates derived therein  show that the method is only of first order convergent in time when SM-ABC  are considered.

In the present work we propose an iterative predictor-corrector method based on the explicit method investigated  in \cite{ABG2017}, resulting a fully explicit method that is second order convergent in time for the SM-ABC case. 
 In the Section \ref{sec:2} we prove that the explicit iterative method converges to a second order in time implicit method and we deduce the {\it a priori} error estimates for the fully discrete scheme.  In Section \ref{sec:3} we illustrate the theoretical results with  some numerical examples and, in Section \ref{sec:4}, we apply the numerical method to a computational model that aims to simulate the light scattering through the outer nuclear layer of the retina.

This work was developed in the framework of a more general project that aims to develop a computational model to simulate the electromagnetic wave's propagation through the eye's structures in order to create a virtual optical coherence  tomography scan \cite{santos2015}.


\section{An explicit iterative leap-frog discontinuous Galerkin method}
 \label{sec:2}
 
 %
%
%
\subsection{Numerical scheme}
\label{scheme_subsec}
 Assume that the computational domain $\Omega$ is a bounded polygonal set that is partitioned into $K$  triangular elements  $T_k$ such that $\overline \Omega = \cup_k T_k$. For simplicity, we consider that the resulting mesh $\mathcal T_h$ is conforming.  
The finite element space is then taken to be $V_{N}=\{v \in L^2(\Omega)^3: v|_{T_k} \in P_N(T_k)^3\},$ 
where $P_N(T_k)$ denotes the space of polynomials of degree less than or equal to $N$ on $T_k$. %
On each element $T_k$, the solution fields $ E_{x}(x,y,t), E_{y}(x,y,t), H_{z}(x,y,t)$ are approximated by the piecewise polynomial functions  
$ E_{xk}(x,y,t),\hat E_{yk}(x,y,t), \hat H_{zk}(x,y,t).$  The approximate fields are allowed to be discontinuous across element boundaries.   In this way, we introduce the notation for the jumps of the field values across the interfaces of the elements, $[\hat E]=\hat E^--\hat E^+$ and
$[\hat H]=\hat H^--\hat H^+$, where the superscript $``+"$ denotes the neighboring element and the superscript $``-"$ refers to the local cell.  Furthermore we introduce, respectively, the cell-impedances and cell-conductances
$Z^{\pm}=\mu^{\pm}c^{\pm}$ and $Y^{\pm}=\left(Z^{\pm}\right)^{-1}$. At the outer cell boundaries we set $Z^+=Z^-$. The coupling between elements is introduced via the numerical flux as in \cite{ABG2017}.

To define the fully discrete scheme, we divide the time interval  into $M$ subintervals by the points $0=t^0<t^1<\cdots<t^M=T$, where $t^m=m \Delta t$, $\Delta t$ is the time step size and $T+ \Delta t/2 \leq T_f$.
The unknowns related to the electric field are approximated at integer time-stations $t^m$ and are denoted by  $\hat E_k^m=\hat E_k(.,t^m)$. The unknowns related to the magnetic field are approximated at half-integer time-stations $t^{m+1/2}=(m+~\frac{1}{2}) \Delta t$ and are denoted by  $\hat H_k^{m+1/2}=\hat H_k(.,t^{m+1/2})$. 

With the above setting , we can now formulate the iterative leap-frog DG method. The process starts with an approximation to the initial data which we denote  by  $(\hat E_{x}^{0}, \hat E_{y}^{0}, \hat H_{z}^{1/2}) \in V_N$.    For each $m=0,1,\ldots, M-1$,  we initialize the iterative process by 
\[\hat E_{x}^{m+1,0}=\hat E_{x}^{m}, \quad \hat E_{y}^{m+1,0}=\hat E_{y}^{m}, \quad \hat H_{z}^{m+3/2,0}=\hat H_{z}^{m+1/2}.\]
The  $(n+1)^{\text{th}}$ inner iteration of  the iterative scheme,   for $n=0,1,2, \ldots$,   is: 
 find $(\hat E_{x}^{m+1,n+1}, \hat E_{y}^{m+1,n+1}, \hat H_{z}^{m+3/2,n+1}) \in V_N$ such that, for all $(u, v, w) \in V_N$ we have 
\begin{eqnarray}
&& \left(\epsilon_{xx} \frac{ \hat E_{x_k}^{m+1,n+1}-\hat E_{x_k}^{m}}{\Delta t} +\epsilon_{xy}  \frac{ \hat E_{y_k}^{m+1,n+1}-\hat E_{y_k}^{m}}{\Delta t} ,u_k \right)_{T_k} = \left(\partial_y \hat H_{z_k}^{m+1/2},
   u_k\right) _{T_k} \nonumber \\
&&\quad+   \left(\frac{-n_y}{Z^++Z^-}\left(Z^+ [\hat H_{z}^{m+1/2}]-\alpha \left(n_x[\hat E_{y}^{[m+1/2,n]}]-n_y[\hat E_{x}^{[m+1/2,n]}]\right)\right) ,u_k\right)_{\partial T_k},  \nonumber \\
\label{Fdis1_it} \\
&&\left(\epsilon_{yx} \frac{ \hat E_{x_k}^{m+1,n+1}-\hat E_{x_k}^{m}}{\Delta t} + \epsilon_{yy} \frac{ \hat E_{y_k}^{m+1,n+1}-\hat E_{y_k}^{m}}{\Delta t}, v_k \right) _{T_k } = -\left( \partial_x \hat H_{z_k}^{m+1/2} , v_k \right)_{T_k} \nonumber \\
&&\quad+ \left(\frac{n_x}{Z^++Z^-} \left(Z^+ [\hat H_{z}^{m+1/2}]-\alpha \left(n_x[\hat E_{y}^{[m+1/2,n]}]-n_y[\hat E_{x}^{[m+1/2,n]}]\right)\right) , v_k\right)_{\partial T_k}, \nonumber \\
\label{Fdis2_it}  \\
&&\left ( \mu \frac{ \hat H_{z_k}^{m+3/2,n+1}-\hat H_{z_k}^{m+1/2}}{\Delta t}, w_k\right)_{T_k} =  \left( \partial_y \hat E_{x_k}^{m+1}
- \partial_x \hat E_{y_k}^{m+1} ,w_k \right)_{T_k} \nonumber \\
&&\quad+ \left (\frac{1}{Y^++Y^-} \left(Y^+ (n_x[\hat E_y^{m+1}]-n_y[\hat E_x^{m+1}])-\alpha[\hat H_z^{[m+1,n]}]\right), w_k\right)_{\partial T_k},  \label{Fdis3_it} 
\end{eqnarray}
where $(\cdot,\cdot)_{T_k}$ and $(\cdot,\cdot)_{\partial T_k}$ denote the classical  $L^2(T_k)$ and $L^2(\partial T_k)$  inner-products and  $\hat {E}^{[m+1/2,n]}$  and $\hat {H}^{[m+1,n]}$  are the average approximations
\[\hat{E}^{[m+1/2,n]}= \frac {\hat{E}^{m}+\hat{E}^{m+1,n}}{2}, \quad \hat{H}^{[m+1,n]}=\frac{\hat{H}^{m+1/2}+\hat{H}^{m+3/2,n}}{2}. \]
The parameter $\alpha \in \{0,1\}$  can be used to control dissipation. Taking $\alpha=0$ yields a non dissipative central flux while $\alpha=1$ corresponds to the classic upwind flux. 

The boundary conditions are discretised  as in \cite{ABG2017, alvarez2014}, this is, for both upwind and central fluxes, consider $\alpha=1$ for the numerical flux at the outer boundary and
$[\tilde E_x]=\tilde E_x^-$, $[\tilde E_y]=\tilde E_y^-$ and $[\tilde H_z]=\tilde H_z^-$.

The current time step $m+1$ is terminated when the stopping criterion 
\[\|\hat E^{m+1,n+1}-\hat E^{m+1,n}\|_{L^2(\Omega)}<\mbox{tol}, \quad \|\hat H_z^{m+3/2,n+1}-\hat H_z^{m+3/2,n}\|_{L^2(\Omega)}<\mbox{tol},\]
is satisfied for some pre-defined small constant $ \mbox{tol}$. Then the correspondent numerical solution is denoted by $(\hat E_{x_k}^{m+1}, \hat E_{y_k}^{m+1}, \hat H_{z_k}^{m+3/2})$.
If we only perform one iteration ($n=0$) we obtain the explicit method considered in \cite{ABG2017}. If we perform two iterations ($n=0$ and $n=1$) we obtain a  predictor-corrector type method.

%
%
%
\subsection{Convergence result}
\label{conv_subsec}
We will show that, under a suitable stability condition, the solution of the iterative predictor-corrector scheme (\ref{Fdis1_it})--(\ref{Fdis3_it}) converges to the solution of the underlying implicit method. The implicit method is defined as follows: given an initial approximation  $(\tilde E_{x}^{0}, \tilde E_{y}^{0}, \tilde H_{z}^{1/2}) \in V_N$, for each $m=0,1,\ldots, M-1$, we compute  $(\tilde E_{x}^{m+1}, \tilde E_{y}^{m+1}, \tilde H_{z}^{m+3/2}) \in V_N$ such that, for all $(u, v, w) \in V_N$,
\begin{eqnarray}
&& \left(\epsilon_{xx} \frac{ \tilde E_{x_k}^{m+1}-\tilde E_{x_k}^{m}}{\Delta t} +\epsilon_{xy}  \frac{ \tilde E_{y_k}^{m+1}-\tilde E_{y_k}^{m}}{\Delta t} ,u_k \right)_{T_k} = \left(\partial_y \tilde H_{z_k}^{m+1/2},
   u_k\right) _{T_k} \nonumber \\
&&\quad+   \left(\frac{-n_y}{Z^++Z^-}\left(Z^+ [\tilde H_{z}^{m+1/2}]-\alpha \left(n_x[\tilde E_{y}^{[m+1/2]}]-n_y[\tilde E_{x}^{[m+1/2]}]\right)\right) ,u_k\right)_{\partial T_k},
 \nonumber \\ \label{Fdis1_imp}\\
&&\left(\epsilon_{yx} \frac{ \tilde E_{x_k}^{m+1}-\tilde E_{x_k}^{m}}{\Delta t} + \epsilon_{yy} \frac{ \tilde E_{y_k}^{m+1}-\tilde E_{y_k}^{m}}{\Delta t}, v_k \right) _{T_k } = -\left( \partial_x \tilde H_{z_k}^{m+1/2} , v_k \right)_{T_k} \nonumber \\
&&\quad+ \left(\frac{n_x}{Z^++Z^-} \left(Z^+ [\tilde H_{z}^{m+1/2}]-\alpha \left(n_x[\tilde E_{y}^{[m+1/2]}]-n_y[\tilde E_{x}^{[m+1/2]}]\right)\right) , v_k\right)_{\partial T_k},
 \nonumber \\ \label{Fdis2_imp}\\
&&\left ( \mu \frac{ \tilde H_{z_k}^{m+3/2}-\tilde H_{z_k}^{m+1/2}}{\Delta t}, w_k\right)_{T_k} =  \left( \partial_y \tilde E_{x_k}^{m+1}
- \partial_x \tilde E_{y_k}^{m+1} ,w_k \right)_{T_k} \nonumber \\
&&\quad+ \left (\frac{1}{Y^++Y^-} \left(Y^+ (n_x[\tilde E_y^{m+1}]-n_y[\tilde E_x^{m+1}])-\alpha[\tilde H_z^{[m+1]}]\right), w_k\right)_{\partial T_k},
\label{Fdis3_imp}
\end{eqnarray}
where we consider   the 
 average approximations  $\tilde{E}^{[m+1/2]}$  for $\tilde{E}^{m+1/2}$ and $\tilde{H}^{[m+1]}$  for $\tilde{H}^{m+1}$ given by 
  \begin{equation}\label{average_approx}
  \tilde{E}^{[m+1/2]}=\frac{\tilde{E}^{m}+\tilde{E}^{m+1}}{2}, \quad  \tilde{H}^{[m+1]}=\frac{\tilde{H}^{m+1/2}+\tilde{H}^{m+3/2}}{2}. 
  \end{equation}
We note that the the numerical solutions are defined  implicitly, since the upwind fluxes involve the unknowns $\tilde E_{x}^{m+1}$, $\tilde E_{y}^{m+1}$ and $\tilde H_{z}^{m+3/2}$.

 Let $h_k$ be the diameter of the triangle $T_k \in \mathcal T_h$, and $h$ be the maximum element diameter, that is, $h_k=\sup_{P_1,P_2 \in T_K} \|P_1-P_2\|$, $h=\max_{T_k \in \mathcal T_h}\{h_k\}.$ 
We assume that the mesh is regular in the sense that  there exists a constant $\tau >0$ such that for all $T_k \in \mathcal T_h$, $\frac{h_k}{\tau_k} \leq \tau$,
where $\tau_k$ denotes the maximum diameter of a ball inscribed in $T_k$. It may be proved (see \cite{riviere2008}) that, for any  $u \in P_N(T_k)$, the following trace inequality holds 
 \begin{equation}\label{inq_trace}
 \|u\|_{L^2(f_k)} \leq C_{\tau} \sqrt{(N+1)(N+2)}  h_k^{-1/2}  \|u\|_{L^2(T_k)},
 \end{equation}
where $f_k$ is   an edge of $T_k$ and $C_{\tau}$ a positive constant independent of $h_k$ and $N$ but dependent on the shape-regularity $\tau$.

 Let us now define the difference between two successive numeric values of the electromagnetic fields by 
\begin{align*}
\delta_n \hat E_{x_k}^{m+1}  &=\hat E_{x_k}^{m+1,n+1}-\hat E_{x_k}^{m+1,n},\\
\delta_n \hat E_{y_k}^{m+1}  &=\hat E_{y_k}^{m+1,n+1}-\hat E_{y_k}^{m+1,n},\\
\delta_n \hat H_{z_k}^{m+3/2} & =\hat H_{z_k}^{m+3/2,n+1}-\hat H_{z_k}^{m+3/2,n}, 
\end{align*}
for $n = 0, 1, 2, \cdots $. 
The following theorem gives upper bounds for $\delta_n \hat E_{x_k}^{m+1}$, $\delta_n \hat E_{y_k}^{m+1}$ and $\delta_n \hat H_{z_k}^{m+3/2}$. 

\begin{theorem}
\label{itv_theorem}
The the solution of the iterative predictor-corrector scheme (\ref{Fdis1_it})--(\ref{Fdis3_it}) converges to the solution of the method (\ref{Fdis1_imp})--(\ref{Fdis3_imp})   provided that the stability condition of the underlying explicit method (i.e.,(\ref{Fdis1_it})--(\ref{Fdis3_it}) taking only the iteration $n=0$)  is satisfied, that is (see \cite{ABG2017})
  \begin{equation}\label{stab_cond_exp}
\Delta t < \frac{\min\{\munderbar \epsilon, \munderbar\mu\}}{\max\{{C}_E, C_H\}} \min\{h_k\},
\end{equation}
  with 
\begin{align*}{C}_E & = \frac{1}{2}C_{inv}N^2+C_{\tau}^2(N+1)(N+2)\left(\frac{5}{2}+\frac{\alpha+\frac{1}{4}}{\min\{Z_k\}}\right),\\
{C}_H & = \frac{1}{2}C_{inv}N^2+ C_{\tau}^2(N+1)(N+2)\left(\frac{5}{2}+\frac{\alpha+\frac{1}{2}}{\min\{ Y_k\}}\right), 
\end{align*}
where $C_{\tau}$ satisfies the trace inequality (\ref{inq_trace}), $C_{inv}$ is a positive constant  independent of $h_k$ and  $N$, and $Z_k$ and $Y_k$ denote respectively the cell-impedance $Z$ and the cell-conductance $Y$ inside the triangle $T_k \in \mathcal T_h$.
\end{theorem}
{\it Proof:} 
The stability condition (\ref{stab_cond_exp}) ensures 
 that $\|\delta_{0}\hat E^{m+1}\|_{L^2(\Omega)}$ and $\|\delta_{0}\hat H_{z}^{m+3/2}\|_{L^2(\Omega)}$ are bounded for all $m=0,1,...,M-1$.

Let us denote by $F^{int}$ the set of internal edges and $F^{ext}$ the set of edges that belong to the boundary $\delta \Omega$. Let $v_k$ be the set of indices of the neighbouring elements of $T_k$. For each $i\in v_k$, we consider the internal edge $f_{ik} = T_i\cap T_k$, and we denote by $n_{ik}$ the unit normal oriented from $T_i$ towards $T_k$. For each boundary edge $f_k = T_k \cap \delta \Omega$, $n_k$ is taken to be the unitary outer normal vector to $f_k$. 

Taking the difference of (\ref{Fdis1_it})--(\ref{Fdis3_it})  between two successive iterations, $n+1$ and $n$, and replacing $u_k$, $v_k$ and $w_k$ by, respectively, $\delta_n \hat E_{x_k}^{m+1}$, $\delta_n \hat E_{y_k}^{m+1}$ and $\delta_n \hat H_{z_k}^{m+3/2}$ and summing over all elements  $T_k \in \mathcal T_h$, we obtain
\begin{align*}
\sum_{T_k \in \mathcal T_h} & \left(\epsilon \delta_n \hat E_{k}^{m+1}  , \delta_n \hat E_{k}^{m+1} \right)_{T_k}   \\
 = &
\frac{ \Delta t}{2}\sum_{f_{ik} \in F^{int} } \int_{f_{ik}} \Bigg(   \frac{(n_y)_{ki}}{Z_i+Z_k}  \left((n_x)_{ki}\delta_{n-1}[\hat E_{y_k}^{m+1}]-(n_y)_{ki}\delta_{n-1}[\hat E_{x_k}^{m+1}]\right) \delta_n \hat E_{x_k}^{m+1}\nonumber\\
 &+ \frac{(n_y)_{ik}}{Z_i+Z_k} \left((n_x)_{ik}\delta_{n-1}[\hat E_{yi}^{m+1}]-(n_y)_{ik}\delta_{n-1}[\hat E_{xi}^{m+1}] \right)
\delta_n \hat E_{x_i}^{m+1}\Bigg)\, ds\nonumber\\
& -\frac{ \Delta t}{2}\sum_{f_{ik} \in F^{int} } \int_{f_{ik}} \Bigg(   \frac{(n_x)_{ki}}{Z_i+Z_k}  \left((n_x)_{ki}\delta_{n-1}[\hat E_{y_k}^{m+1}]-(n_y)_{ki}\delta_{n-1}[\hat E_{x_k}^{m+1}]\right)
\delta_n \hat E_{y_k}^{m+1}\nonumber\\
&  +\frac{(n_x)_{ik}}{Z_i+Z_k} \left((n_x)_{ik}\delta_{n-1}[\hat E_{yi}^{m+1}]-(n_y)_{ik}\delta_{n-1}[\hat E_{xi}^{m+1}]\right)
\delta_n \hat E_{y_i}^{m+1}\Bigg)\, ds \nonumber\\
&+\frac{ \Delta t}{2}\sum_{f_{k} \in F^{ext} } \int_{f_{k}}\Bigg(\frac{(n_y)_{k}}{2 Z_k}\left((n_x)_{k}\delta_{n-1}[\hat E_{y_k}^{m+1}]-(n_y)_{k}\delta_{n-1}[\hat E_{x_k}^{m+1}]\right)
\delta_{n}\hat E_{x_k}^{m+1}\nonumber\\
 & -\frac{(n_x)_{k}}{2 Z_k}  \left((n_x)_{k}\delta_{n-1}[\hat E_{y_k}^{m+1}]-(n_y)_{k}\delta_{n-1}[\hat E_{x_k}^{m+1}]\right)
\delta_{n}\hat E_{y_k}^{m+1}\Bigg)\, ds,
\end{align*}
\begin{align*}
 \sum_{T_k \in \mathcal T_h}  & \left(  \mu \delta_{n} \hat H_{z_k}^{m+3/2},\delta_{n} \hat H_{z_k}^{m+3/2}\right)_{T_k}  =  \\ 
 & -\frac{ \Delta t}{2}\sum_{f_{ik} \in F^{int} } \int_{f_{ik}} \Bigg(\frac{1}{Y_i+Y_k} \delta_{n-1}[\hat H_{z_k}^{m+3/2}]\delta_{n}\hat H_{z_k}^{m+3/2}
+\frac{1}{Y_i+Y_k} \delta_{n-1}[\hat H_{zi}^{m+3/2}]\delta_{n}\hat H_{zi}^{m+3/2}
\Bigg)\, ds\\
& -\frac{ \Delta t}{2}\sum_{f_{k} \in F^{ext} } \int_{f_{k}}\Bigg( \frac{1}{2 Y_k} \delta_{n-1}[\hat H_{z_k}^{m+3/2}]
\delta_{n}\hat H_{z_k}^{m+3/2}\Bigg)\, ds.
\end{align*}
Then
\begin{align*}
\sum_{T_k \in \mathcal T_h} & \left(\epsilon \delta_n \hat E_{k}^{m+1}  , \delta_n \hat E_{k}^{m+1} \right)_{T_k} \\ 
= &
\frac{ \Delta t}{2}\sum_{f_{ik} \in F^{int} } \int_{f_{ik}} \Bigg(   \frac{(n_y)_{ki}}{Z_i+Z_k}  \left((n_x)_{ki}\delta_{n-1}[\hat E_{y_k}^{m+1}]-(n_y)_{ki}\delta_{n-1}[\hat E_{x_k}^{m+1}]\right) \delta_n [\hat E_{x_k}^{m+1}]\Bigg)\, ds\nonumber\\
& - \frac{ \Delta t}{2}\sum_{f_{ik} \in F^{int} } \int_{f_{ik}} \Bigg(   \frac{(n_x)_{ki}}{Z_i+Z_k}  \left((n_x)_{ki}\delta_{n-1}[\hat E_{y_k}^{m+1}]-(n_y)_{ki}\delta_{n-1}[\hat E_{x_k}^{m+1}]\right)
\delta_n [\hat E_{y_k}^{m+1}]\Bigg)\, ds \nonumber\\
& +\frac{ \Delta t}{2}\sum_{f_{k} \in F^{ext} } \int_{f_{k}}\Bigg(\frac{(n_y)_{k}}{2 Z_k}\left((n_x)_{k}\delta_{n-1}[\hat E_{y_k}^{m+1}]-(n_y)_{k}\delta_{n-1}[\hat E_{x_k}^{m+1}]\right)
\delta_{n}\hat E_{x_k}^{m+1}\nonumber\\
 &  \qquad \qquad \qquad \quad -\frac{(n_x)_{k}}{2 Z_k}  \left((n_x)_{k}\delta_{n-1}[\hat E_{y_k}^{m+1}]-(n_y)_{k}\delta_{n-1}[\hat E_{x_k}^{m+1}]\right),
\delta_{n}\hat E_{y_k}^{m+1}\Bigg)\, ds, 
\end{align*}
and 
\begin{align*}
\sum_{T_k \in \mathcal T_h} & \left( \mu \delta_{n} \hat H_{z_k}^{m+3/2},\delta_{n} \hat H_{z_k}^{m+3/2}\right)_{T_k} \\ = &
-\frac{ \Delta t}{2}\sum_{f_{ik} \in F^{int} } \int_{f_{ik}} \Bigg(\frac{1}{Y_i+Y_k} \delta_{n-1}[\hat H_{z_k}^{m+3/2}]\delta_{n}[\hat H_{z_k}^{m+3/2}]\Bigg)\, ds\\
& -\frac{ \Delta t}{2}\sum_{f_{k} \in F^{ext} } \int_{f_{k}}\Bigg( \frac{1}{2 Y_k} \delta_{n-1}[\hat H_{z_k}^{m+3/2}]
\delta_{n}\hat H_{z_k}^{m+3/2}\Bigg)\, ds.
\end{align*}
So
\begin{align*}
\sum_{T_k \in \mathcal T_h} & \left(\epsilon \delta_n \hat E_{k}^{m+1}  , \delta_n \hat E_{k}^{m+1} \right)_{T_k} \nonumber\\  & \leq 
\frac{ \Delta t}{4\min\{Z_k\}}\sum_{f_{ik} \in F^{int} } \|\delta_{n-1}[\hat E_{k}^{m+1}]\|_{L^2(f_{ik})} \|\delta_n [\hat E_{k}^{m+1}]\|_{L^2(f_{ik})}\nonumber\\
& \quad +\frac{ \Delta t}{4\min\{Z_k\}}\sum_{f_{k} \in F^{ext} }
\|\delta_{n-1}[\hat E_{k}^{m+1}]\|_{L^2(f_{k})}
\|\delta_n [\hat E_{k}^{m+1}]\|_{L^2(f_{k})},
\end{align*}
\begin{align*}
\sum_{T_k \in \mathcal T_h} & \left( \mu \delta_{n} \hat H_{z_k}^{m+3/2},\delta_{n} \hat H_{z_k}^{m+3/2}\right)_{T_k} \nonumber\\ &  \leq
\frac{ \Delta t}{4\min\{Y_k\}}\sum_{f_{ik} \in F^{int} } \|\delta_{n-1}[\hat H_{z_k}^{m+3/2}]\|_{L^2(f_{ik})} \|\delta_{n}[\hat H_{z_k}^{m+3/2}]\|_{L^2(f_{ik})}\\
&\quad  +\frac{ \Delta t}{4\min\{Y_k\}}\sum_{f_{k} \in F^{ext} } \|\delta_{n-1}[\hat H_{z_k}^{m+3/2}]\|_{L^2(f_{k})}
\|\delta_{n}\hat H_{z_k}^{m+3/2}\|_{L^2(f_{k})}.
\end{align*}
Consequently, considering (\ref{inq_trace}), we obtain
\begin{align*}
&\left(\munderbar\epsilon-\frac{ \Delta t}{\min\{Z_k\}} C_{\tau}^2 (N+1)(N+2)\max\left\{h_k^{-1} \right\}\right) \| \delta_n \hat E^{m+1} \|_{L^2(\Omega)} \\
&\qquad \leq
 \frac{ \Delta t}{\min\{Z_k\}} C_{\tau}^2 (N+1)(N+2)\max\left\{h_k^{-1} \right\} \|\delta_{n-1}\hat E^{m+1}\|_{L^2(\Omega)},
\end{align*}
\begin{align*}
&\left(\munderbar\mu-\frac{ \Delta t}{\min\{Y_k\}} C_{\tau}^2 (N+1)(N+2)\max\left\{h_k^{-1} \right\}\right) \| \delta_{n} \hat H_{z_k}^{m+3/2}\|_{L^2(\Omega)}\\ &\qquad \leq
\frac{ \Delta t}{\min\{Y_k\}} C_{\tau}^2 (N+1)(N+2)\max\left\{h_k^{-1} \right\} \|\delta_{n-1}\hat H_{z}^{m+3/2}\|_{L^2(\Omega)}.
\end{align*}

Taking the following condition into account (that results from (\ref{stab_cond_exp}))
\[
\Delta t < \frac{\min\{\munderbar \epsilon, \munderbar\mu\}\min\{Z_k, Y_k\}}{C_{\tau}^2 (N+1)(N+2)} \min\{h_k\},
\]
we conclude the proof. $\Box$

The next theorem establishes that the implicit method is second order convergent in time and arbitrary high order in space and so, with the previous result, we may conclude that the same occurs for the iterative explicit scheme.

 \begin{theorem} \label{theo_convergence}
Let us consider the implicit leap-frog DG method (\ref{Fdis1_imp})--(\ref{Fdis3_imp}) complemented with the discrete boundary conditions defined in Section \ref{scheme_subsec} and suppose that the solution of the Maxwell's equations (\ref{MaxwellTE}) complemented by  (\ref{SM}) has the following regularity:
$E_x, E_y, H_z \in L^\infty(0,T_f;H^{s+1}(\Omega))$,  $\frac{\partial E_x}{\partial t}, \frac{\partial E_y}{\partial t}, \frac{\partial H_z}{\partial t} \in L^2 (0,T_f;H^{s+1}(\Omega) \cap L^\infty(\partial\Omega))$ and $\frac{\partial^2 E_x}{\partial t^2},$  $\frac{\partial^2 E_y}{\partial t^2},$ $\frac{\partial^2 H_z}{\partial t^2}$ $\in$ 
$L^2(0,T_f;H^{1}(\Omega))$, $s \geq 0$. If  the time step $\Delta t$ satisfies
\begin{equation}\label{stab_cond_imp}
\Delta t <  \frac{\min\{\munderbar \epsilon, \munderbar\mu\}}{\frac12 C_{inv}N^2 + 2C_{\tau}^2(N+1)(N+2)} \min\{h_k\} ,
\end{equation}
where $C_{inv}$ and $C_{\tau}$ are the positive constants defined in the previous theorem, then
\begin{eqnarray*}
\max_{1\leq m\leq M} \left(\|E^{m}-\tilde E^m\|_{L^2(\Omega)}+\|H_{z}^{m+1/2}-\tilde H_{z}^{m+1/2}\|_{L^2(\Omega)} \right)
\leq C ( \Delta t^2 + h^{\min\{s,N\}})\\
+ C  \left(\|E^{0}-\tilde E^0\|_{L^2(\Omega)}+\|H_{z}^{1/2}-\tilde H_{z}^{1/2}\|_{L^2(\Omega)} \right)
\end{eqnarray*}
holds, where $C$ is a generic positive constant independent of $\Delta t$ and the mesh size~$h$.
\end{theorem}
 
 {\it Proof:}  Follows the steps of the proof of Theorem 4.2 in \cite{ABG2017}. $\Box$


\section{Numerical results}
\label{sec:3}

To  illustrate the theoretical results of the previous section, we consider the model problem (\ref{MaxwellTE}) defined in the  square $\Omega = (-1,1)^2$, complemented by initial conditions and SM-ABC (\ref{SM}).  In order to  make it  easier to find examples with known
 exact solution and consequently with the possibility to compute the error of the numerical solution,  source terms were introduced in the model.
  The simulation time is fixed at $T = 1$ and in all tests we set $\mu = 1$ and $\epsilon$ is given by (\ref{tensor}), with $\epsilon_{xx}=4x^2+y^2+1$, $\epsilon_{yy}=x^2+1$ and $\epsilon_{xy}=\epsilon_{yx}=\sqrt{x^2+y^2}$. The source terms are defined in such way that the problem has the exact solution
\begin{align*}E_x(x,y,t)  & = -\sqrt{\frac{\epsilon_{yy}}{\det(\epsilon)}}\sin(\pi t)\sin(\pi x),\\
 E_y(x,y,t)  & =  \sqrt{\frac{\epsilon_{xx}}{\det(\epsilon)}}\sin(\pi t)\sin(\pi y), \\
 H_z(x,y,t ) &  =  \sin(\pi t) \sin(\pi xy).
 \end{align*}

\begin{figure}[ht]
\includegraphics[scale=.27]{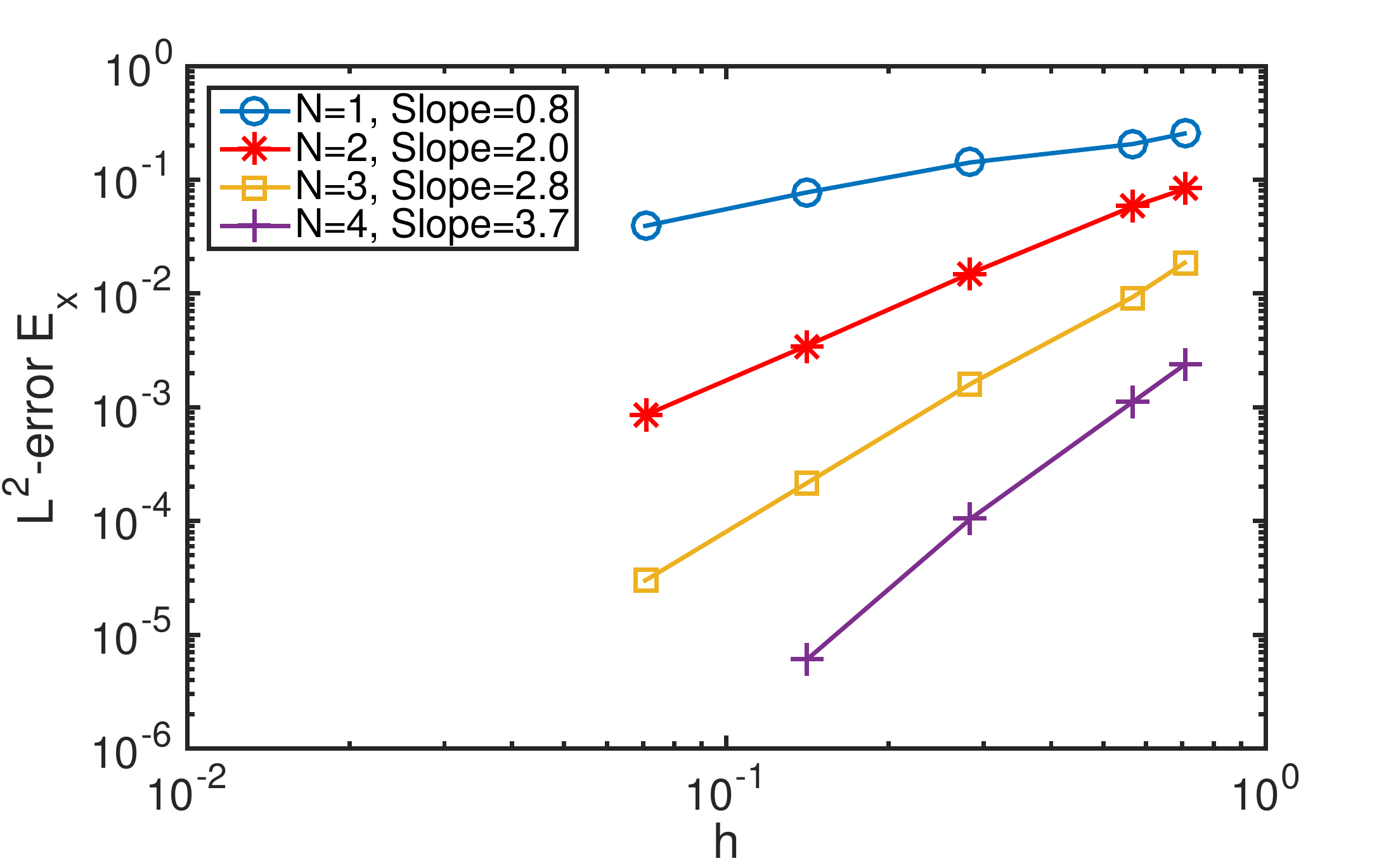}\ 
\includegraphics[scale=.27]{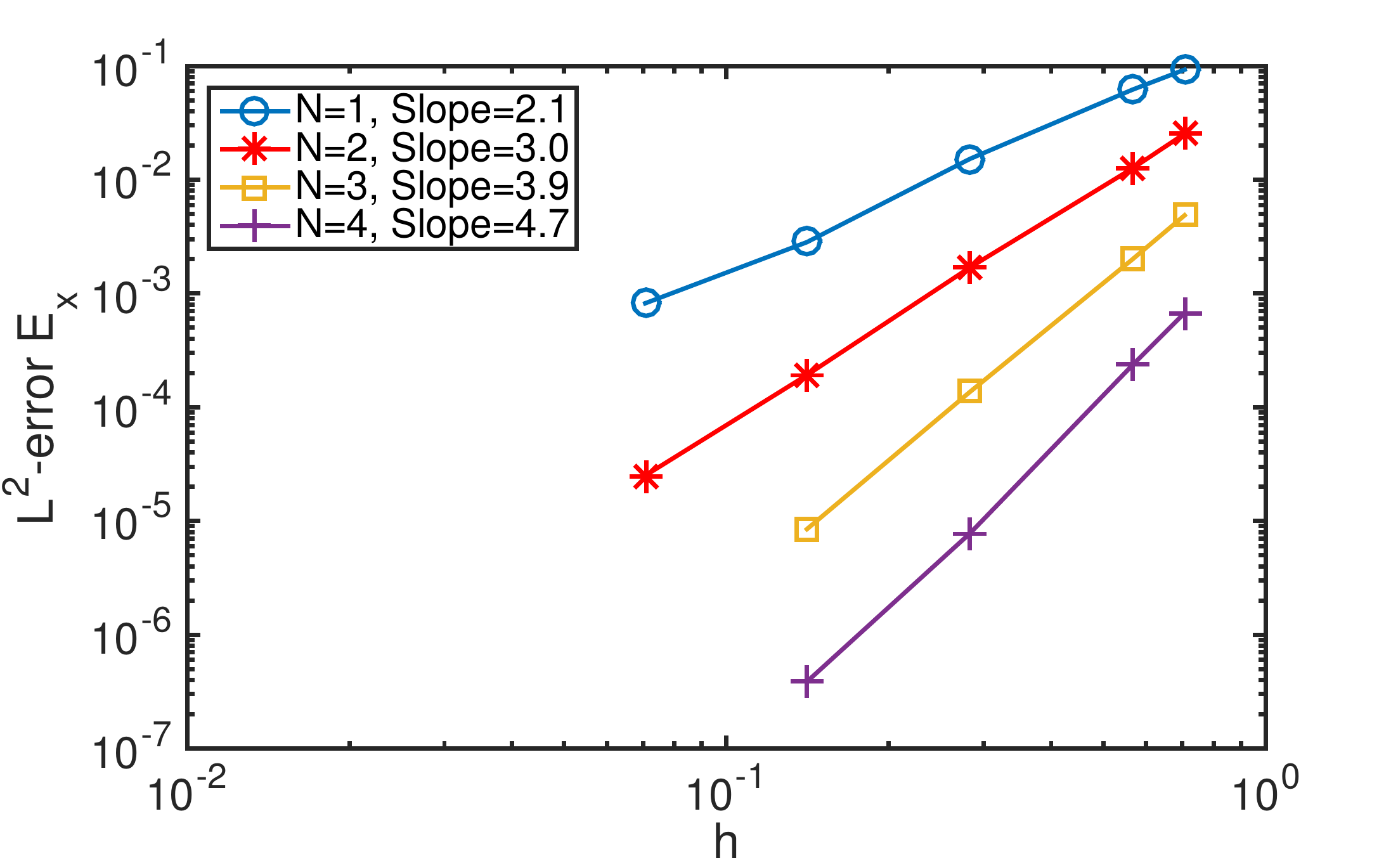}
\caption{$\|E_x^{M}-\tilde E_x^M\|_{L^2(\Omega)}$ {\it versus } $h$. Left: central flux; Right: upwind flux}
\label{fig:space}      
\end{figure}

  To illustrate the order of convergence in space, we fix $\Delta t = 10^{-5 }$, except when $N=4$ where we consider $\Delta t=10^{-6}$. 
In Fig. \ref{fig:space} we plot the discrete $L^2$-error of the $\tilde E_x$ component of electric field depending on the maximum element diameter of each mesh,  for different degrees for the polynomial approximation, for both central and upwind fluxes. The vertical and horizontal axis are scaled logarithmically. The numerical convergence rate is approximated by the slope of the linear regression line. 
    For central flux, the numerical convergence rate is close to the value estimated  in Theorem \ref{theo_convergence},  $\mathcal{O}(h^{N})$,  and for upwind flux we observe higher order  of convergence, up to  $\mathcal{O}(h^{N+1})$ in some cases.
Similar results were obtained for $\tilde E_y$ and $\tilde H_z$.

\begin{figure}[ht]
\includegraphics[scale=.27]{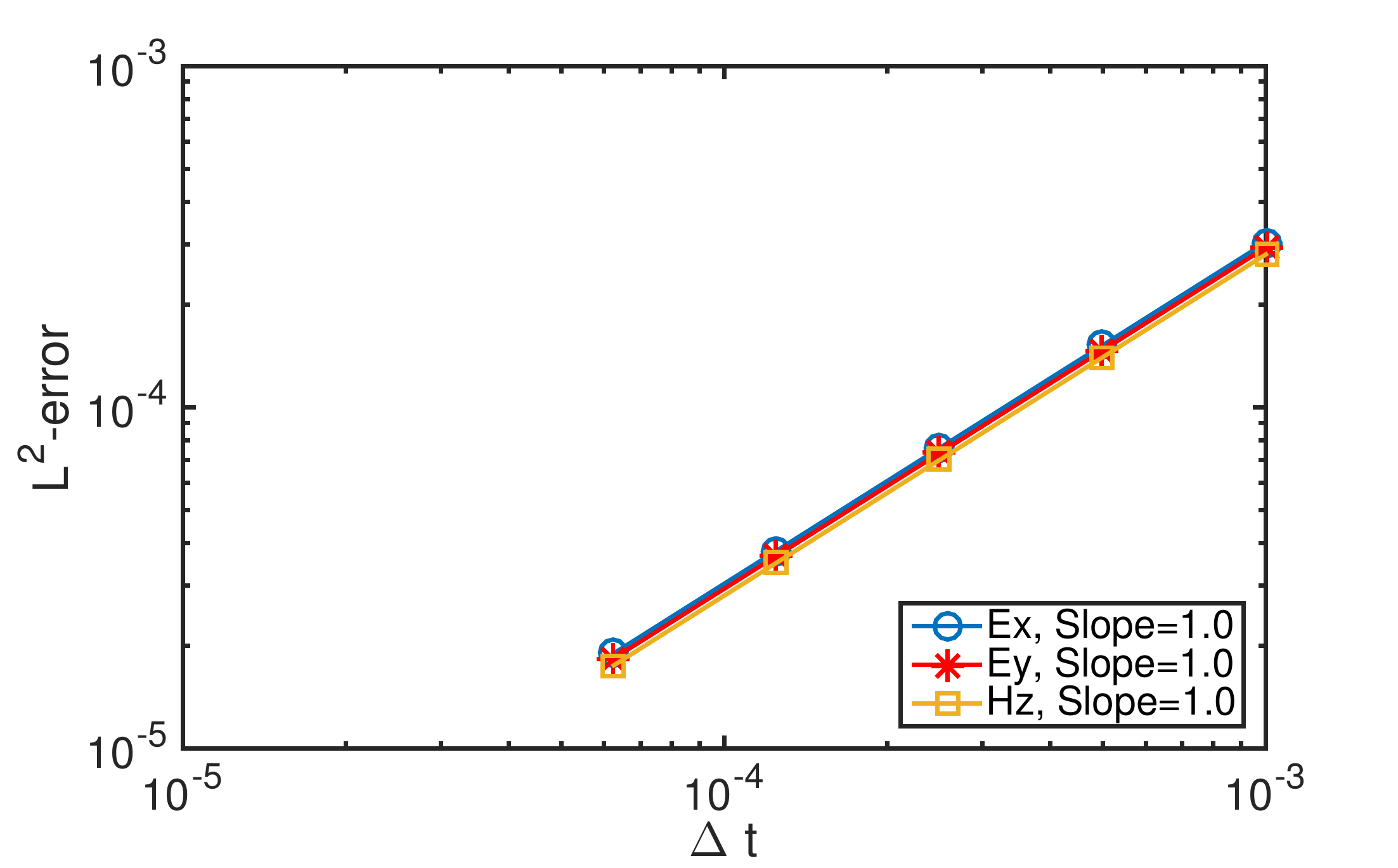}\ 
\includegraphics[scale=.27]{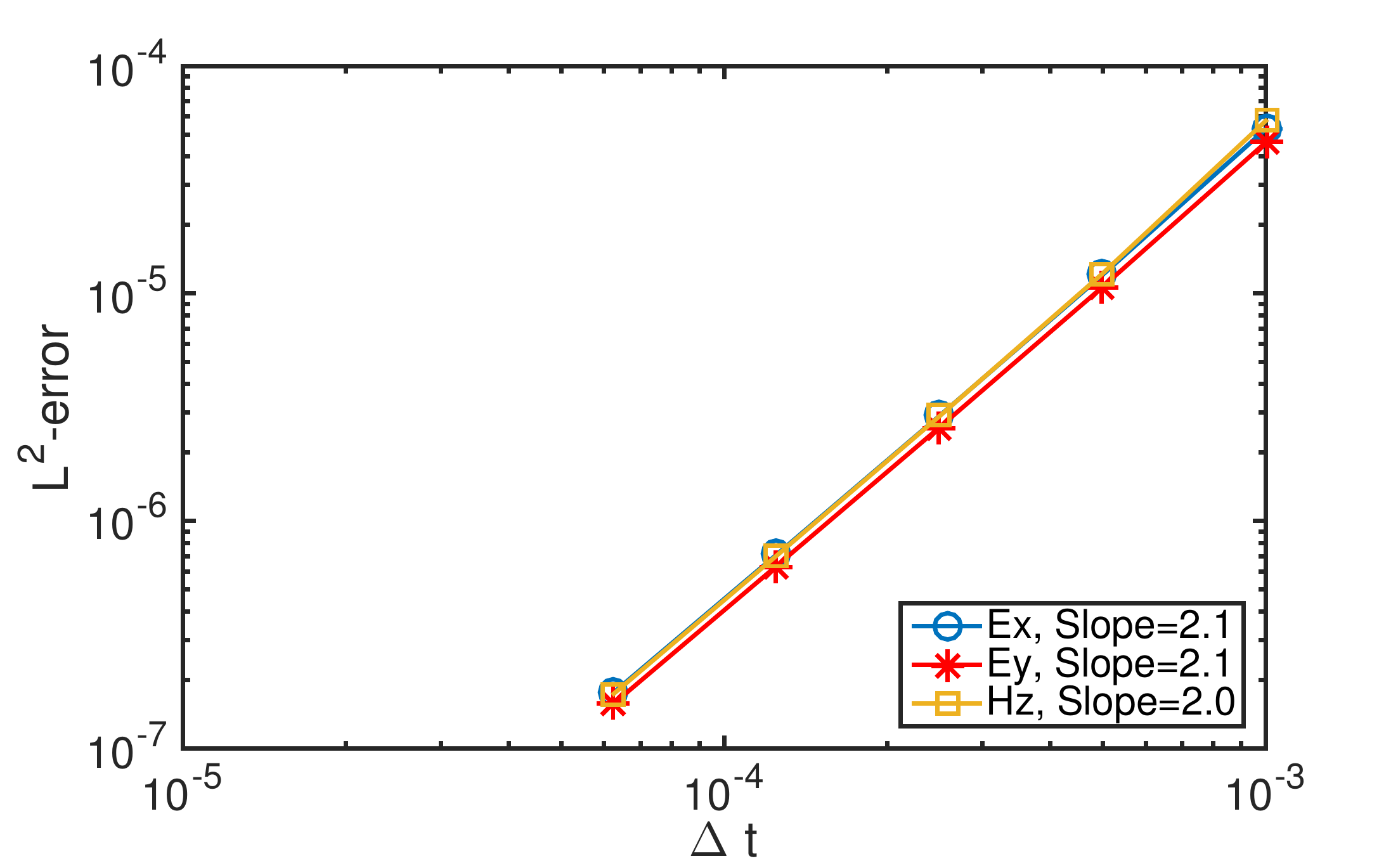}
\caption{$\|E_x^{M}-\tilde E_x^M\|_{L^2(\Omega)}$, $\|E_y^{M}-\tilde E_y^M\|_{L^2(\Omega)}$  and $\|H_z^{M}-\tilde H_z^M\|_{L^2(\Omega)}$ {\it versus } $\Delta t$. Left: explicit leap-frog method; Right: predictor-corrector  leap-frog method.}
\label{fig:time}   
\end{figure}

To visualize the convergence order in time, the  polynomial degree and the number of elements have been  set to $N = 8$ and $K = 800$, respectively. The results plotted  in Fig. \ref{fig:time} illustrate the  first order  convergency in time for the explicit leap-frog method and  show that the second order is recovered when the predictor-corrector method is considered. These results correspond to upwind fluxes. The experiments using central fluxes show analogous results in terms of order of convergence in time. 


\section{Modeling scattered electromagnetic wave's propagation through eye's structures}
\label{sec:4}

This work is part of a research project which aims to develop a cellular model of the human retina able to simulate different retinal/cellular conditions and how these changes are translated to an Optical Coherence Tomography scan  \cite{santos2015}. Simulating the full complexity of the retina, in particular the variation of the size and shape of each structure, distance between them and the respective refractive indexes, requires a rigorous approach that can be achieved by solving Maxwell's equations. As the interest is to acquire the backscattered light intensity, we start this section by the scattered  field formulation. Then we build up a two dimensional model which tries to represent a single nucleus of the outer nuclear layer (ONL) of the retina. The performance of our method is examined by simulating the light scattering in this 2D domain. The evolution of the scattering  field intensity in time is obtained using the predictor-corrector DG method.

\subsection{Optical Coherence Tomography (OCT)}

\begin{figure}[h!]
\begin{center}
\includegraphics[scale=.5]{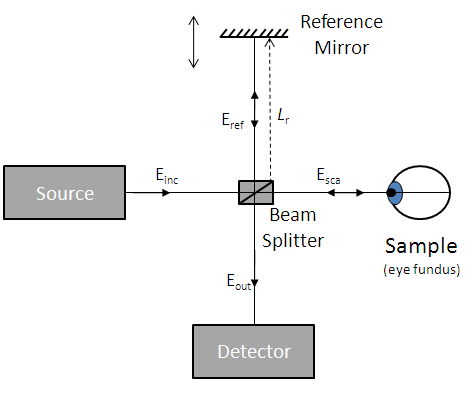}
\end{center}
\caption{Scheme for the principle of OCT \cite{serranho2012}.}
\label{oct}   
\end{figure}

Optical Coherence Tomography (OCT) its an imaging technology that produces high- resolution cross-sectional images of the internal microstructure of living tissue, widely used in ophthalmological exams. This technology's working principle is analogous to ultrasound, but it uses light instead of sound to locate subtle differences in the tissue being analysed. Discontinuities in the refractive index of the tissue give rise to light scattering, with some light backscattered to the detector. Factors such as the shape and size of the scatterer, wavelength of the incident light and refractive index differences have an impact on the amount of backscattered light. During a scan, the OCT machine directs a light beam into the retina and extracts, through interferometry, the backscattered light intensity of retinal structures and their depth location in an A-scan (see Fig. \ref{oct}). By transversely moving the light beam, several A-scans can be collected into a cross-sectional image -- a B-scan. Usually, several cross- sectional images are acquired by probing an azimuthal direction and combined into a volume.

\subsection{The scattered field formulation}
We can exploit the linearity of the Maxwell's equations  in order to separate the electromagnetic  fields ($E$, $H$) into incident fields  ($E^i$, $H^i$) and scattered components ($E^s$, $H^s$), {\it{i.e.}}, $E = E^s + E^i$ and $H = H^s + H^i$.

Assuming that the incident field is also a solution of the Maxwell's equations we obtain in the same way as in \cite{taflove1995},  the scattered field formulation,
\begin{align}
  \epsilon_{xx}\frac{\partial E_x^{s}}{\partial t}  +  \epsilon_{xy}\frac{\partial E_y^{s}}{\partial t} &= \frac{\partial H_z^{s}}{\partial y} + P,\label{Maxwell SF1}\\
  \epsilon_{yx}\frac{\partial E_x^{s}}{\partial t}   + \epsilon_{yy}\frac{\partial E_y^{s}}{\partial t} &= -\frac{\partial H_z^{s}}{\partial x} + Q,\label{Maxwell SF2}\\
  \mu\frac{\partial H_z^{s}}{\partial t} &= -\frac{\partial E_y^{s}}{\partial x} + \frac{\partial E_x^{s}}{\partial y} + R,\label{Maxwell SF3} \qquad \textrm{in }  \quad \Omega\times(0,T], 
\end{align}
with the source terms
\begin{align}
P(x,y,t) &= (\epsilon^i - \epsilon_{xx})\frac{\partial E_x^{i}}{\partial t}  - \epsilon_{xy}\frac{\partial E_y^{i}}{\partial t},\\
Q(x,y,t) &= - \epsilon_{yx}\frac{\partial E_x^{i}}{\partial t} + (\epsilon^i - \epsilon_{yy})\frac{\partial E_y^{i}}{\partial t},\\
R(x,y,t) &= (\mu^i - \mu)\frac{\partial H_z^{i}}{\partial t},
\end{align}
where $\epsilon^i$ and $\mu^i$ represent, respectively,  the relative permittivity and permeability of the medium in which the incident field propagates in the absence of scatterers (in the background medium).   
Additionally, using this formulation it is straightforward to specify an incident wave using an analytic formula.

The intensity of the light  that hits the OCT detector defines the output signal. Hence, we are interested in computing 
the scattered field intensity,   
\begin{equation}
\label{intensity_function}
I^s=\sqrt{(E_x^s)^2+(E_y^s)^2}.
\end{equation} 


\subsection{Light scattering  in the outer nuclear layer}
 We use our numerical model to simulate light scattering in the ONL. This layer has a special relevance among the retina's layers as it consistently presents the characteristics of  diabetic macular edema  \cite{ciulla2005, correia2014}.  The ONL is mostly populated by the cells bodies of light sensitive photoreceptor cells (rods and cons). Thus, we postulate that the main contribution to light scattering in this layer comes from the nucleus \cite{seet2009}, as it is the biggest organelle in the soma and presents a high refractive index difference to the surrounding medium. As such, the ONL can be modelled  as a population of spherical nuclei in an homogenous medium. 
As a proof of concept we present a simple simulation in a two  dimensional square domain which contains circles that aims to represent, respectively, a single nucleus and three nuclei in the ONL. The permittivity inside the circles and in the background domain has different values.

Let us consider equations (\ref{Maxwell SF1})--(\ref{Maxwell SF3}), in $\Omega=(-1,1)^2$, complemented with SM-ABC  and null initial conditions.
The absorbing boundary conditions are chosen for the model as they avoid undesirable reflections that invade the computational domain. 
In the first experiment we will consider the case of just one circle:  
${\cal{C}} = \{(x,y)\in \Omega : x^2+y^2<0.25\}.$
In the second example we will consider the case of three circles: 
${\cal{C}}_1 = \{(x,y)\in \Omega : x^2+(y-0.5)^2<0.01\};$
${\cal{C}}_2 = \{(x,y)\in \Omega :  x^2+y^2<0.01\};$
${\cal{C}}_3 = \{(x,y)\in \Omega : x^2+(y+0.5)^2<0.01\}.$
In the experiments the relative permittivity and permeability  and magnetic permeability are  considered as constants, $\epsilon^i=1$ and  $\mu^i=\mu=1$. 
The electric permittivity  is considered as  a diagonal matrix with $\epsilon_{xx}(x,y)=\epsilon_{yy}(x,y)=1.2$ for $(x,y)$ inside the circles and  $\epsilon_{xx}(x,y)=\epsilon_{yy}(x,y)=1$ otherwise. 
For the  incident wave we consider the planar wave $E_y^i(x,t) = \cos(10(x-t))$. 
\begin{figure}[h!]
\begin{center}
\includegraphics[scale=.33]{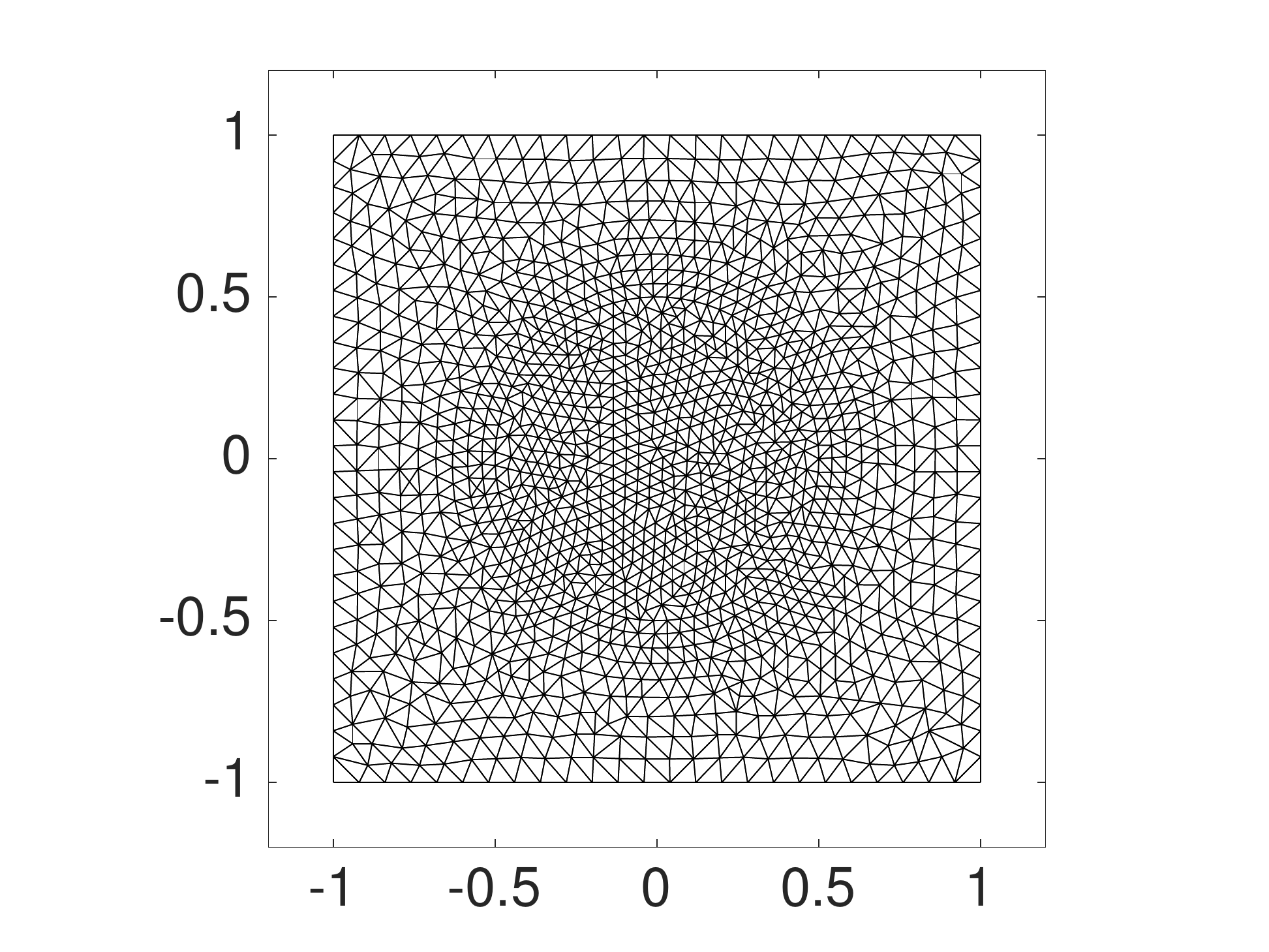}\hspace{-3.5em}
\includegraphics[scale=.33]{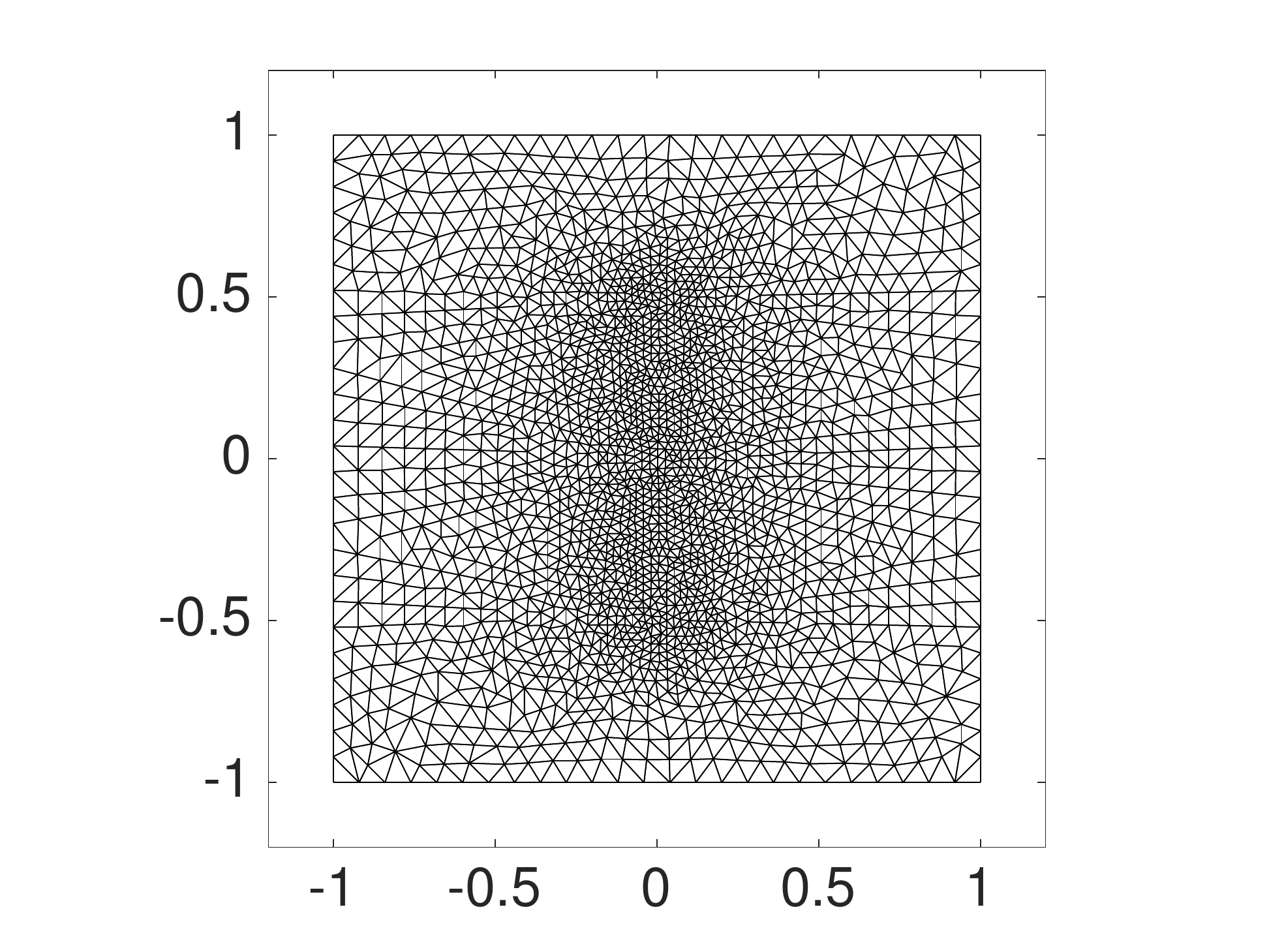}
\end{center}
\caption{Square computational domain $\bar \Omega =[-1,1]^2$  and the  triangular mesh used in the computations: one circle (left); three circles (right).}
\label{mesh_circle}   
\end{figure}
%

%
%
 \begin{figure}[h!]
\includegraphics[width=0.45\textwidth]{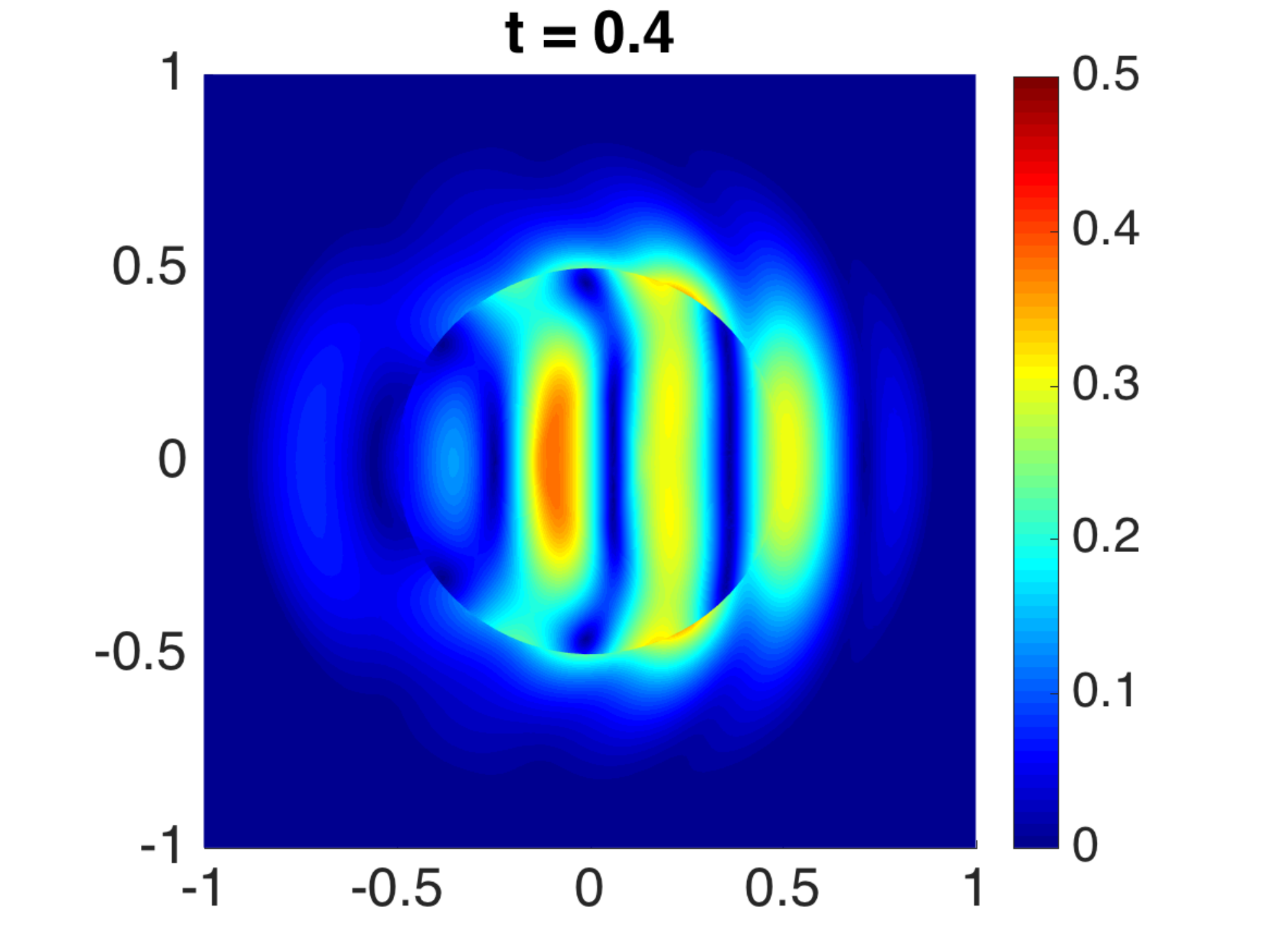}\hspace{-1em}
\includegraphics[width=0.45\textwidth]{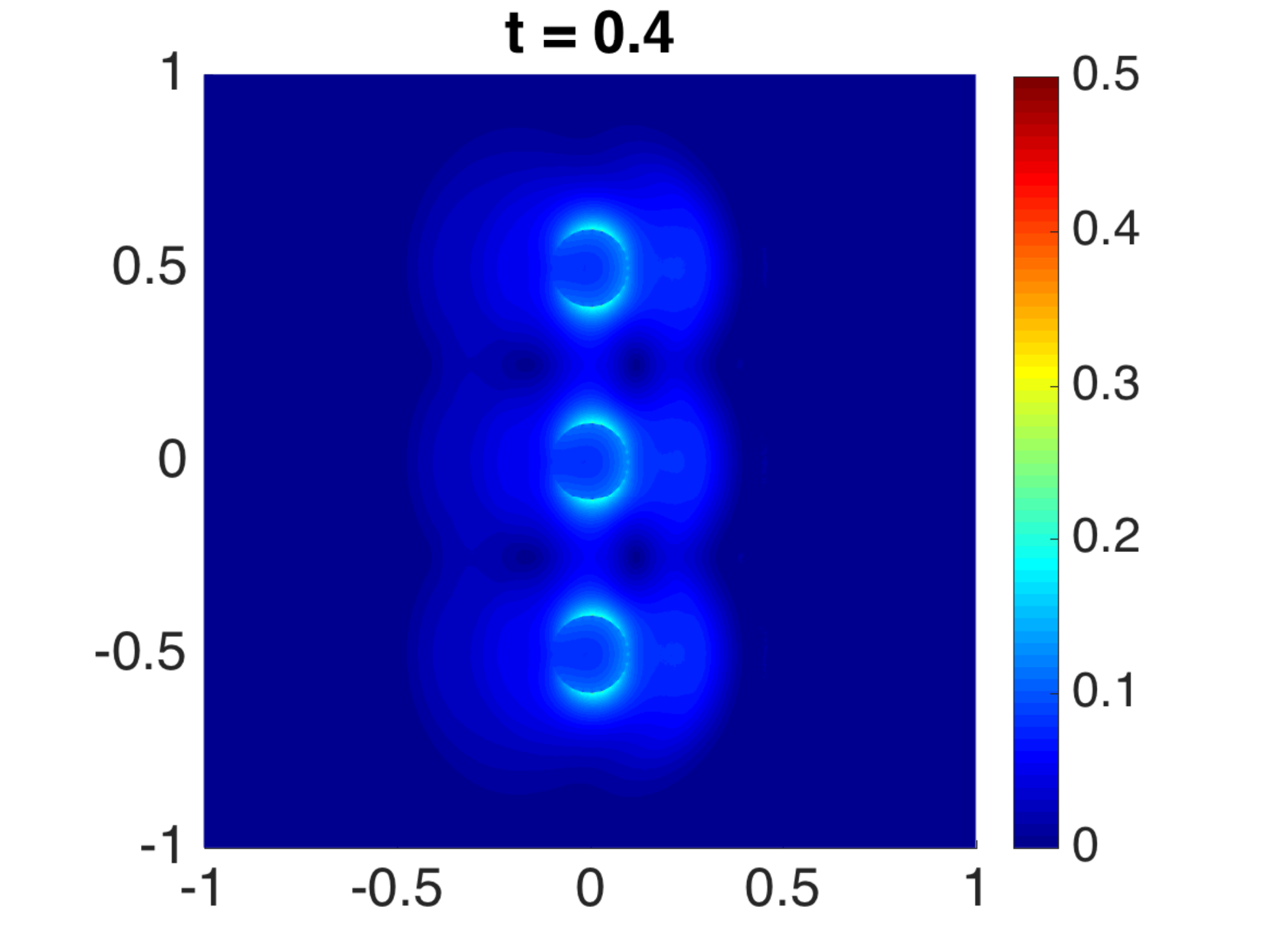}
\includegraphics[width=0.45\textwidth]{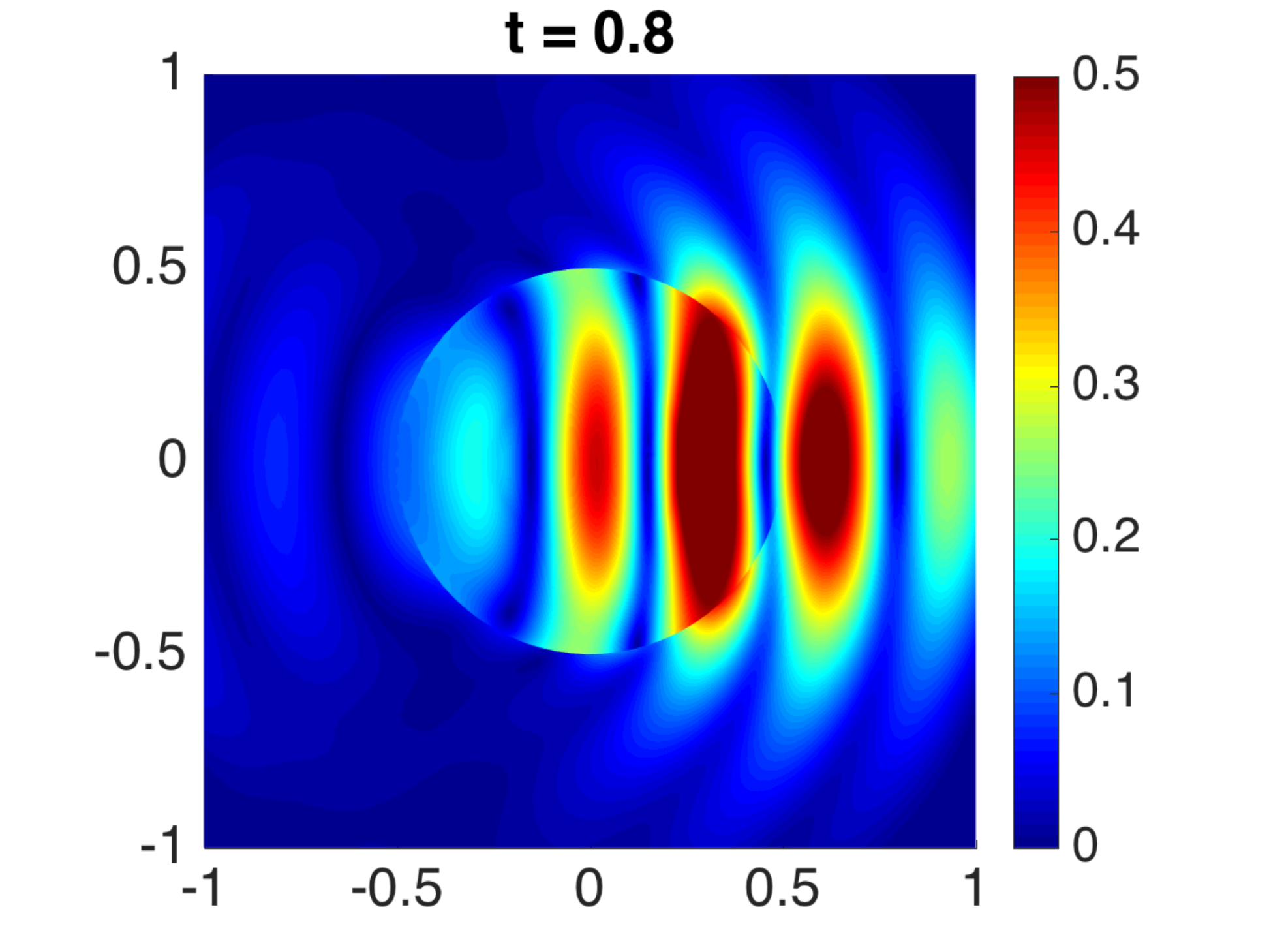}\hspace{-1em}
\includegraphics[width=0.45\textwidth]{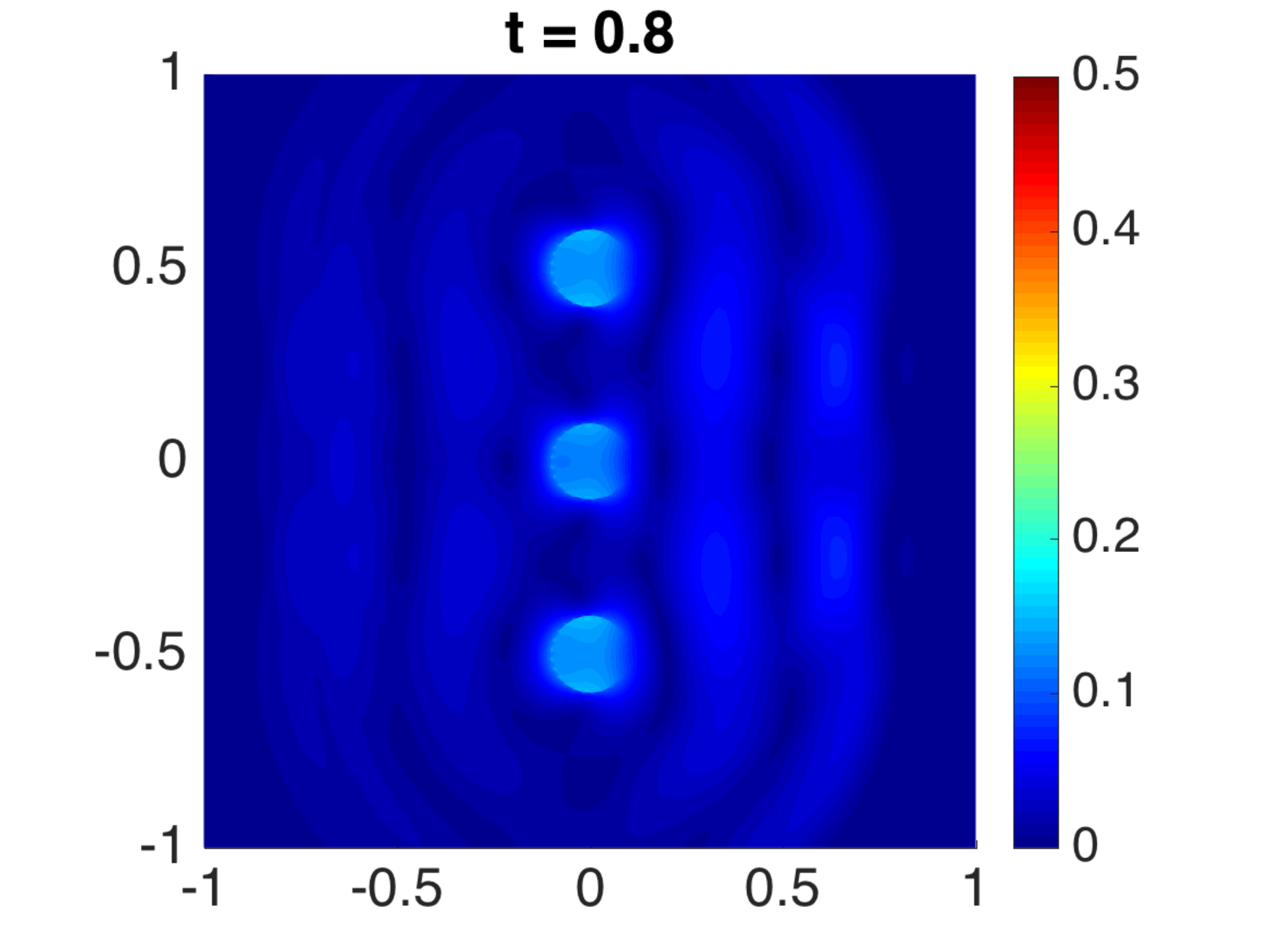}
\caption{Evolution of the scattered field intensity (\ref{intensity_function}) with time: one circle (left); three circles (right).}
\label{intensity}
\end{figure}

 For the simulation we used with predictor-corrector DG method defined in Section \ref{sec:2}, considering $\alpha=0$ (central flux) and the approximation polynomial degree    $N=4$. The time step was chosen to be  $\Delta t=0.002$ and the final simulation time is $T = 0.8$.  The meshes are illustrated in Fig.~\ref{mesh_circle}. 
The evolution in time of the  scattered field intensity (\ref{intensity_function}) 
 is plotted in  Fig.~\ref{intensity}. These results show that the scatterers are  clearly identified. With this model, we can simulate more complex cellular structures only by changing  the electric permittivity tensor $\epsilon$.

\section{Conclusions}
\label{sec:5}

We presented an iterative explicit  leap-frog DG method for time dependent Maxwell's equations in anisotropic media, considering SM-ABC.  The numerical scheme is fully explicit and converges to a second order in time implicit method. 
The results of a set of numerical experiments  support  the  theoretical results.  Moreover we developed a 2D model which simulates the light scattering by  a single nucleus in the outer nuclear layer of the retina. This work was elaborated in the framework of a more general project with a real application (see \cite{correia2014, santos2015}).


\begin{backmatter}

\section*{Competing interests}
  The authors declare that they have no competing interests.

\section*{Author's contributions}
All authors contributed equally to the writing of this paper. All authors read and approved the final manuscript.

\section*{Acknowledgements}

This work was partially supported by the Centre for Mathematics of the University of Coimbra -- UID/MAT/00324/2013, funded by the Portuguese Government through FCT/MCTES and co-funded by the European Regional Development Fund through the Partnership Agreement PT2020; by the Portuguese Government through  the BD grant SFRH/BD/51860/2012; and by the grant LARSyS  UID/EEA/50009/2013.

\end{backmatter}
\end{document}